\documentclass{ifacconf}
\usepackage{amsmath}
\usepackage{mathrsfs,amssymb,amsmath,booktabs,array,xcolor,url,cite,color,soul,multirow,multicol}
\usepackage{mathbbold,bbm}
\usepackage{stfloats}
\usepackage{amsfonts}
\usepackage{extarrows}
\usepackage{textcomp}
\usepackage{ulem}
\usepackage{cases}
\usepackage{bm}
\usepackage{arydshln}
\usepackage[all]{xy}
\usepackage{ulem}
\usepackage{graphicx}      
\usepackage{natbib}
\usepackage{extarrows}
\usepackage{widetext}
\usepackage{flushend}
\usepackage{cuted}
\usepackage{graphicx}
\usepackage{epstopdf}
\newtheorem{theorem}{Theorem}
\newtheorem{definition}{Definition}
\newtheorem{remark}{Remark}

\newtheorem{proposition}{Proposition}
\begin{document}
\begin{frontmatter}

\title{Lyapunov Function Partial Differential Equations for Stability Analysis of a Class of Chemical Reaction Networks}

\thanks[footnoteinfo]{This work was supported by
the National Natural Science Foundation of China [Grant No. 11671418,
11271326 and 61611130124].}

\author[First]{Shan~Wu}
\author[First]{Yafei~Lu}
\author[First]{Chuanhou~Gao }

\address[First]{School of Mathematical Sciences, Zhejiang University, Hangzhou 310027, China (e-mail: 11335036@zju.edu.cn (S. Wu); 11535029@zju.edu.cn (Y. Lu); gaochou@zju.edu.cn (C. Gao))}

\begin{abstract}                
We investigate a broad family of chemical reaction networks (CRNs) assigned with mass action kinetics, called complex-balanced-produced-CRNs (CBP-CRNs), which are generated by any given complex balanced mass action system (MAS)
and whose structures depend on the selection of producing matrices. Unluckily, the generally applied pseudo-Helmholtz free energy function may fail to act as a Lyapunov function for the CBP-CRNs. Inspired by the method of Lyapunov function partial differential equations (PDEs), we construct one solution of their corresponding Lyapunov function PDEs, termed as the generalized pseudo-Helmholtz free energy function, and further show that solution can behave as a Lyapunov function to render the asymptotic stability for the CBP-CRNs. This work can be taken as an argument of the conjecture that Lyapunov function PDEs approach can serve for any MAS.
\end{abstract}

\begin{keyword}
CBP-CRNs, Lyapunov function PDEs, the generalized pseudo-Helmholtz free energy function
\end{keyword}

\end{frontmatter}

\section{Introduction}
Chemical reaction networks (CRNs) play unique roles in chemistry, biology, process industry and many other fields. The dynamics of CRNs following the mass action law, known as the mass action systems (MASs), can be modeled as polynomial ordinary differential equations. Plenty of work has been done on characterizing the dynamical properties of such systems [\cite{Angeli2007A}, \cite{Angeli2015A}, \cite{Feinberg1972Complex}, \cite{Finberg1987Chemical}, \cite{HornJackson1972General}], among which the Lyapunov stability analysis has drawn great attention [\cite{Alradhawi2016New}, \cite{Anderson2015Lyapunov}, \cite{HornJackson1972General}, \cite{Rao2013A}, \cite{Sontag2012Structure}]. For certain MASs with special structures such as reversible, balanced, etc., some inspiring outcomes have been achieved [\cite{Anderson2010Product}, \cite{Craciun2006Multiple}, \cite{Feinberg1995The}, \cite{Johnston2013Computing}, \cite{Rao2013A}]. Following these studies, this paper deals with the construction of Lyapunov functions for MASs with some special structures.\\
It has been showed that some particular network structures of MASs can lead to proper properties of equilibrim stability and distribution. Considering the weakly reversible structure, \cite{HornJackson1972General} have discovered that with a zero deficiency such an MAS is complex balanced and admits exactly one positive equilibrium for each positive stoichiometric compatibility class, meanwhile each of these equilibria is locally asymptotically stable, known as the deficiency zero theorem.  \cite{Finberg1988Chemical} has then further worked out the deficiency one theorem, weakened the deficiency conditions for the former one. Moreover, given persistent property of a complex balanced MAS, the global stability of its equilibria is proved [\cite{Anderson2011A},\cite{Craciun2013Persistence},  \cite{Pantea2011On}, \cite{Siegel2000Global}, \cite{Sontag2012Structure}]. Such studies have naturally drawn attention to the special structure of complex balancing.\\
For Lyapunov stability analysis of CRNs, the construction of proper Lyapunov functions is a critical step. Sorts of approaches have been applied on this study, including employing the scaling limits of nonequilibrium potential basing on the potential theory[\cite{Anderson2015Lyapunov}], piecewise linear in rate Lyapunov functions [\cite{Alradhawi2016New}] for balanced MASs, etc.. Specified to complex balanced MASs, \cite{HornJackson1972General} have found the pseudo-Helmohotlz free energy function as a proper candidate. Another way related to the complex balancing is working on the linear conjugacy [\cite{Johnston2011Linear}, \cite{Johnston2012Dynamical}, \cite{Johnston2013Computing}], where it has been showed that linear conjugate systems share stability properties, thus a complex balanced conjugate of the considered MAS proves its locally asymptotical stability [\cite{Szederke2011Finding}]. \cite{Ke2019Complex} developed the generalized pseudo-Helmholtz function as the Lyapunov function for some conjugate MAS defined as reconstructions and reverse reconstructions of CRNs. Meanwhile, the Lyapunov functions for CRNs with more general structures remain challenging. In this case, \cite{Fang2015Lyapunov} derived a partial differential equation (PDE) based on chemical master equation for any general MAS, with its solutions (if existed) as Lyapunov candidates of such system, and thus referred to as the Lyapunov function PDEs. This method has served well on complex balanced CRNs, general CRNs with a 1-dimensional stoichiometric subspace and some special CRNs with higher dimensional stoichiometric subspace, while general CRNs with more than a 2-dimensional stoichiometric subspace still need to be explored. \cite{Fang2015Lyapunov} have made a conjecture that for any MAS admitting a stable positive equilibrium, the induced Lyapunov function PDEs have the required solution under proper set of boundary complex. Inspired by these works, this paper focuses on a broad family of CRNs generated by certain complex balanced MASs, where common Lyapunov candidates including the pseudo-Helmholtz function fail for the stability analysis. We develop the Lyapunov function PDEs for such systems, and seek for their proper solutions to serve as Lyapunov functions. Such CRNs also inherit some neat properties on asympotically stability and equilibria distribution from complex balanced networks. This work extends the Lyapunov stability analysis to a larger range of MASs.\\
The remainder of this paper is organized as follows. Section 2 prepares some basic concepts on CRNs and Lyapunov function PDEs. Section 3 starts with some conclusions on applying Lyapunov function PDEs to complex balanced MASs, then by linear conjugacy approach, from complex balanced MASs we generate a class of networks named complex-balanced-produced-CRNs (CBP-CRNs), for which the Lyapunov function PDEs can be solved by some function with similar forms of the pseduo-Helomholtz function. This function is proved to serve well as a Lyapunov candidate for the CBP-CRNs. After this, Section 4 illustrate our results by an example. Finally, Section 5 concludes this work.


~\\~\\
\noindent{\textbf{Mathematical Notations:}}\\
\rule[1ex]{\columnwidth}{0.8pt}
\begin{description}
\item[\hspace{-0.5em}{$\mathbb{R}^n:$}] $n$-dimensional real space.
\item[\hspace{-0.5em}{$ \mathbb{R}^n_{\geq 0}:$}] $n$-dimensional non-negative real space.
\item[\hspace{-0.5em}{$\mathbb{R}^n_{>0}:$}] $n$-dimensional positive real space.
\item[\hspace{-0.5em}{$\mathbb{Z}^n_{\geq 0}$}]: $n$-dimensional non-negative integer space.
\item[\hspace{-0.5em}{$x^{v_{\cdot i}}$}]: $x^{v_{\cdot i}}=\prod_{j=1}^{d}x_{j}^{v_{ji}}$, where $x, v_{\cdot i}\in \mathbb{R}^{d}$ and $0^{0}=1$.
\item[\hspace{-0.5em}{$\frac{x}{y}$}]: $\frac{x}{y}=(\frac{x_1}{y_1}, \cdots, \frac{x_n}{y_n}) $, where $x\in\mathbb{R}^{n}$, $y\in\mathbb{R}^{n}_{>0}$.
 \item[\hspace{-0.5em}{$\mathrm{Ln}(x)$}]: $\mathrm{Ln}(x)=\left(\ln{x_{1}}, \cdots, \ln{{x}_{n}} \right)^{\top}$, where $x\in\mathbb{R}^{n}_{>{0}}$.
\end{description}
\rule[1ex]{\columnwidth}{0.8pt}
\\~
\section{Basic Concept}
In this section, some basic concepts of CRNs and Lyapunov function PDEs are recalled.
\subsection{chemical reaction networks}
Consider a network with $n$ species $S_1, \cdots, S_n$ and $r$ chemical reactions $\mathcal{R}_1, \cdots, \mathcal{R}_r $ with the ith reaction written as
\begin{equation*}
\mathcal{R}_i:~~~\sum^{n}_{j=1}v_{ji}S_j \rightarrow
\sum^{n}_{j=1}v'_{ji}S_j,
\end{equation*}

where $v_{ji}, ~v'_{ji}\in\mathbb{Z}_{\geq 0}$ represent the complexes of reactant and resultant, then it is defined as follows [\cite{Feinberg1995The}].
\begin{definition}\label{df:basic}
Let $\mathcal{S}=\{S_1, \cdots, S_n\}$, $\mathcal{C}=\bigcup^r_{i=1}\{v_{\cdot i},v_{\cdot i}^{'}\}$ , and $\mathcal{R}=\{v_{\cdot 1}\rightarrow v_{\cdot 1}^{'},\cdots,v_{\cdot r}\rightarrow v_{\cdot r}^{'}\}$, representing the sets of species, complexes, and reactions respectively. When satisfying
\begin{enumerate}
    \item{$\forall v_{\cdot i}\rightarrow v_{\cdot i}^{'}\in\mathcal{R},v_{\cdot i}\neq v_{\cdot i}^{'}$;}
    \item{$\forall S_{j}\in \mathcal{S},v_{\cdot i}\in\mathcal{C}$,the  $j$th entry of $v_{\cdot i}$ represents the stoichiometric coefficient of species $S_{j}$ in complex $v_{\cdot i}$}
\end{enumerate}
the triple $(\mathcal{S,C,R})$ is a \emph{chemical reaction network (CRN)}. \\
Moreover, if the CRN follows the mass action law, i.e, given the concentration vector $x\in\mathbb{R}_{\geq 0}^{n}$ and the rate constant $k_{i}\in\mathbb{R}_{>0}$ for reaction $v_{\cdot i}\rightarrow v_{\cdot i}^{'}$, the reaction rate can be calculated as $k_{i}x^{v_{\cdot i}}$, letting $\mathcal{K}=(k_{1},\cdots,k_{r})$ representing the set of reaction rate constants, then $(\mathcal{S,C,R,K})$ is a \emph{mass action system (MAS)}.
\end{definition}

\begin{definition}
For a CRN $(\mathcal{S,C,R})$, $\mathscr{S}:=span\{v_{\cdot 1}-v_{\cdot 1}^{'},\cdots,v_{\cdot r}-v_{\cdot r}^{'}\}$ is the \emph{stoichiometric subspace} of the network. $\forall C\in\mathbb{R}_{\geq0}^{n}$, $C+\mathscr{S}:=\{C+\xi\mid\xi\in\mathscr{S}\}$ is a \emph{stoichiometric compatibility class} of $C$ for the network.
\end{definition}

Denoting $\Gamma\in\mathbb{Z}_{n\times r}$, $\Gamma _{\cdot i}=v_{\cdot i}^{'}-v_{\cdot i}$ as the stoichiometric matrix, and $R(x)\in \mathbb{R}^{r}$, $R_{i}(x)=k_{i}x^{v_{\cdot i}}$, the dynamics of an MAS$(\mathcal{S,C,R,K})$ can thus be expressed as
\begin{equation}\label{eq:mas}
    \frac{\mathrm{d}x}{\mathrm{d}t}=\Gamma R(x),~~~x\in \mathbb{R}_{\geq0}^{n} ,
\end{equation}
\begin{definition}
For an MAS$(\mathcal{S,C,R,K})$, if there is a concentration $x^{*}\in\mathbb{R}_{>0}^{n}$ s.t. $\Gamma R(x^{*})=0$, then $x^{*}$ is an \emph{equilibrium} of the MAS.\\
If $\exists x^{*}\in\mathbb{R}_{>0}^{n}$ s.t.
\begin{equation}\label{eq:complex}
    \sum_{\{i\mid v_{\cdot i}=z\}} k_{i}(x^{*})^{v_{\cdot i}}= \sum_{\{i\mid v_{\cdot i}^{'}=z\}} k_{i}(x^{*})^{v_{\cdot i}},~~~\forall z\in\mathcal{C}
\end{equation}
which means at certain state the consuming rate equals the producing rate for any complex, then $x^{*}$ is a \emph{complex balanced equilibrium} of the MAS, and $(\mathcal{S,C,R,K})$ is a \emph{complex balanced system}.
\end{definition}
Complex balanced MASs have been proved to possess some elegant properties [\cite{Rao2013A}].
\begin{theorem}\label{thm:complex}
Given any complex balanced MAS ($\mathcal{S}, \mathcal{C}, \mathcal{R}, \mathcal{K}$) with an equilibrium $x^{*}\in \mathbb{R}^{n}_{>0}$, for any initial state $x_{0}\in \mathbb{R}^{n}_{>0}$, there exists a unique positive equilibrium $x^{+}\in({x}_{0}+\tilde{\mathscr{S}})$, and $x^{+}$ is locally asymptotically stable with respect to all initial states in $x_{0}\in \mathbb{R}^{n}_{>0}$ nearby $x^{+}$.
\end{theorem}

\subsection{Lyapunov function PDEs}
From [\cite{Fang2015Lyapunov}], the chemical master equation of a CRN can be developed into a PDE
\begin{equation}\label{eq:LyapunovPDE}
\sum^{r}_{i=1}k_{i}x^{v_{\cdot i}}-\sum^{r}_{i=1}k_{i}x^{v_{\cdot i}}\exp\big\{ (v_{\cdot i}^{'}-v_{\cdot i})^{\top}\nabla f(x)\big\}=0,
\end{equation}
where $x\in\mathbb{R}_{>0}^{n}$, and for which the boundary condition is

\begin{equation}\label{eq:boundary}
\sum_{\{i\mid v_{\cdot i}\in\mathcal{C}_{\bar{x}}\}}k_{i}x^{v_{\cdot i}}-\sum_{\{i\mid v_{\cdot i}^{'}\in\mathcal{C}_{\bar{x}}\}}k_{i}x^{v_{\cdot i}} \exp\big\{(v_{\cdot i}^{'}-v_{\cdot i})^{\top}\nabla f(x)\big\}=0,
\end{equation}
\begin{equation*}
  x\rightarrow \bar{x}, x\in (\bar{x}+\mathscr{S})\cap\mathbb{R}_{>0}^{n}
\end{equation*}
where $\mathcal{C}_{\bar x}$ represents the complex set induced by the boundary point $\bar{x}$.\\

\begin{definition}
For an MAS$(\mathcal{S,C,R,K})$, (\ref{eq:LyapunovPDE}) and (\ref{eq:boundary}) are \emph{Lyapunov function PDEs} for the system.
\end{definition}
One feature of the Lyapunov PDEs is that their solutions (if exist) possess dissipation, that is $\dot {f}(x)\leq 0$ with equality holding if and only if $\nabla f(x) \bot \mathscr{S}$. It implies the potential of the Lyapunov function PDEs to generate Lyapunov function candidates for the MASs with some moderate conditions, as listed below in Theorem \ref{asymstability}. Therefore, the Lyapunov stability analysis can be approached by constructing a proper solution to them.
\begin{theorem}\label{asymstability}
let $x^* \in \mathbb{R}^{n}_{>0}$ be any equilibrium of the MAS ($\mathcal{S}, ~\mathcal{C}, ~\mathcal{R}, ~\mathcal{K}$), if the corresponding Lyapunov function PDEs have a twice differentiable solution $f(x)$ satisfying that $\forall\mu\in\mathscr{S}$, $\mu\nabla^{2}f(x)\mu\geq 0$ with the equality holding if and only if $\mu=\Vec{0}$
all over a region $\mathcal{N}(x^*)=\delta(x^*)\cap(x^*+\mathscr{S})\cap \mathbb{R}_{>0}$ where $\delta ({x^*})$ is a neighborhood of $x^*$, then $f(x)$ can behave as a Lyapunov function to render this system to be locally asymptotically stable at $x^*$ with any initial conditions in $\mathcal{N}(x^*)\cap\{x\mid f(x)<inf_{\{y\in\partial\mathscr{N}(x^*)\}}f(y)\}$.
\end{theorem}

~\\
\section{Main results}
In this section, we aim to acquire the Lyapunov function for a kind of special networks assigned with mass action kinetics, named CBP-CRNs, through the method of Lyapunov function PDEs.
\subsection{Lyapunov function PDEs for complex balanced MASs}
First of all, we will exhibit how Lyapunov function PDEs works for complex balanced MASs. Assume that a complex balanced MAS possesses one equilibrium point, it is well known that this network locally asymptotically converges to the equilibrium taking the pseudo-Helmholtz free energy function as the Lyapunov function [\cite{HornJackson1972General}].
As a matter of fact, the pseudo-Helmholtz free energy function is proved to be one of the solutions of the corresponding Lyapunov function PDEs for complex balanced systems. In addition, this solution serves as a Lyapunov function for the stability analysis for such MASs.
The following Theorem \ref{complexPDEth} is given to indicate the above conclusions.

\begin{theorem}\label{complexPDEth}
For a complex balanced MAS ($\mathcal{S}, \mathcal{C}, \mathcal{R}, \mathcal{K}$) with an equilibrium $x^*\in \mathbb{R}^{n}_{>0}$, its corresponding Lyapunov function PDEs which are described by (\ref{eq:LyapunovPDE}) and (\ref{eq:boundary})
have a solution
\begin{equation}
G(x)=\sum^{n}_{j=1}(x^*_j-x_j-x_j\ln{\frac{x^*_j}{x_j}})
\end{equation}
which is usually called the pseudo-Helmholtz free energy function. This means
\begin{equation}\label{complexPDE}
\sum^{r}_{i=1}k_i x^{v_{\cdot i}}\left(1-\exp \{(v'_{\cdot i}-v_{\cdot i})^\top\nabla G(x)\}\right)=0
\end{equation}
and the boundary condition satisfies \\
\begin{eqnarray}\label{boundcon}
&&\sum_{\{i|v_{\cdot i}\in \mathcal{C}_{\bar x}\}}k_i x^{v_{\cdot i}}
-\sum_{\{i|v'_{\cdot i}\in \mathcal{C}_{\bar x}\}}k_i x^{v_{\cdot i}}\exp \{(v'_{\cdot i}-v_{\cdot i})^\top\nabla G(x)\}\notag \\
&&=0,
 ~~~~~x\rightarrow \bar{x}, ~x\in (\bar{x}+\mathscr{S})\cap\mathbb{R}_{>0}^{n}
\end{eqnarray}

whatever the boundary complex
set ${\mathcal{C}_{\bar{x}}}$ is.
Furthermore, this system can achieve locally asymptotic
stability at $x^*$ for any $x \in (x^*+\mathscr{S})\cap \mathbb{R}^{n}_{>0}$ near $x^*$ with respect to $G(x)$.
\end{theorem}
\begin{pf}
The detailed proof can be found in the literature  [\cite{Fang2015Lyapunov}].
\end{pf}

\subsection{Lyapunov function PDEs for a class of CBP-CRNs}

Referred to [\cite{Fang2015Lyapunov}], the method of Lyapunov function PDEs works not only for complex balanced MASs as above, but also for CRNs with 1-dimensional stoichiometric subspace and some special cases assigned with higher dimension. Still, this approach hasn't been tested on many other structures of CRNs. This paper aims to apply such method to a class of networks which can be generated by complex balanced MASs through linear conjugacy related approach, and are defined as follows.

The CBP-CRNs are defined as
\begin{definition}\label{df:CBP-CRNs}
Given a complex balanced MAS ($\mathcal{S}, \mathcal{C}, \mathcal{R}, \mathcal{K}$) governed by (\ref{eq:mas}) with $x^* \in \mathbb{R}^{n}_{>0}$ to be an equilibrium, a MAS ($\tilde{\mathcal{S}}, \tilde{\mathcal{C}}, \tilde{\mathcal{R}}, \tilde{\mathcal{K}}$) is called a \emph{CBP-CRN} with respect to ($\mathcal{S}, \mathcal{C}, \mathcal{R}, \mathcal{K}$) if for some positive diagonal matrix $D=\text{diag}(d_1,~\cdots, ~d_n)$, its complex set $\tilde{\mathcal{C}}=\cup^{\tilde{r}}_{i=1}\{\tilde{v}_{\cdot i}, ~\tilde{v}'_{\cdot i}\}$
and reaction set $\tilde{\mathcal{R}}=\cup^{\tilde{r}}_{i=1}\{\tilde{v}_{\cdot i}\stackrel{\tilde{k}_i}\longrightarrow\tilde{v}'_{\cdot i}\}$ satisfy
\begin{enumerate}
\item{ $\tilde{r}=r$;}
\item{ $\tilde{v}_{\cdot i}=v_{\cdot i},~ \tilde{v}'_{\cdot i}=v_{\cdot i}+D^{-1}(v'_{\cdot i}-v_{\cdot i}), ~\tilde{k}_i =k_{i}\prod^{n}_{j=1}d^{v_{ji}}_{j}$, $\forall ~i=1, ~\cdots, ~\tilde{r}$.}
\end{enumerate}
\end{definition}

\begin{remark}
It should be pointed out that the constructed CBP-CRNs can be seen as an extension of the networks derived by reverse reconstruction for a given MAS referred to [\cite{Ke2019Complex}] where the original MAS and reverse reconstructed MAS share the same state variables and rate constants.
\end{remark}

\begin{remark}\label{rm:mark1}
Following Definition \ref{df:CBP-CRNs}, it is easy to see that one  complex balanced network can generate a family of corresponding networks with different producing matrices Ds, which indicates that the CBP-CRNs is a class of networks that covers the original network. When D is selected as an unit matrix, the CBP-CRN is exactly the original complex balanced system.
\end{remark}

\begin{remark}
It should be noted that $D$ needs to be well chosen to guarantee $v'_{\cdot i}\in \mathbb{Z}^n_{>0}$ following our MAS setting in Definition \ref{df:basic}. When extending to the generalized mass action systems where $v_{\cdot i}, ~v'_{\cdot i}$ could be fractional, the selection of $D$s can thus be broadened.
\end{remark}

As a straightforward consequence from Definition \ref{df:CBP-CRNs}, we have
\begin{proposition}\label{prop}
For a CBP-CRN ($\tilde{\mathcal{S}}, \tilde{\mathcal{C}}, \tilde{\mathcal{R}}, \tilde{\mathcal{K}}$) generated by a given complex balanced MAS ($\mathcal{S}, \mathcal{C}, \mathcal{R}, \mathcal{K}$) following (\ref{eq:mas}), the dynamical equations of the CBP-CRN follow
\begin{equation}\label{CBP-CRNs}
\dot{\tilde{x}}=\tilde{\Gamma}\tilde{R}(\tilde{x}),
\end{equation}
with $\tilde{\Gamma}=D^{-1}\Gamma$ and $\tilde{R}(\tilde{x})=R(x)$ for $\tilde{x}=D^{-1}x$. Meanwhile, the existence of an equilibrium $x^{*}$ of ($\mathcal{S}, \mathcal{C}, \mathcal{R}, \mathcal{K}$) is equivalent to $\tilde{x}^{*}=D^{-1}x^{*}$ being an equilibrium of ($\tilde{\mathcal{S}}, \tilde{\mathcal{C}}, \tilde{\mathcal{R}}, \tilde{\mathcal{K}}$), and $\tilde{\mathscr{S}}=span\{D^{-1}(v_{\cdot i}-v_{\cdot i}^{'})\mid v_{\cdot i}-v_{\cdot i}^{'}\in\mathscr{S}\}\doteq D^{-1}\mathscr{S}$.
\end{proposition}

Apparently, $G(\tilde{x})$ could not serve as a Lyapunov function anymore for this type of systems described by (\ref{CBP-CRNs}) except the special case where the derived CBP-CRNs gratify the complex balanced condition. In order to obtain a proper Lyapunov function for the CBP-CRN system, we turn our focus to its
corresponding Lyapunov function PDEs, which are written as
\begin{equation}\label{D-PDE}
\sum^{\tilde{r}}_{i=1}\tilde{k}_{i}\tilde{x}^{\tilde{v}_{\cdot i}}\left(1-\exp\{(\tilde{v}'_{\cdot i}-\tilde{v}_{\cdot i})^\top \nabla \tilde{f}(\tilde{x})\}\right)=0,
\end{equation}
together with the boundary condition
\begin{eqnarray}\label{boundcon1}
&&
\sum_{\{i|v_{\cdot i}\in \mathcal{C}_{\hat x}\}}\tilde{k}_i \tilde{x}^{v_{\cdot i}}
-\sum_{\{i|v'_{\cdot i}\in \mathcal{C}_{\hat x}\}}\tilde{k}_i \tilde{x}^{v_{\cdot i}}\exp \{(\tilde{v}'_{\cdot i}-\tilde{v}_{\cdot i})^\top\nabla \tilde{f}(\tilde{x})\}\notag \\
&&=0,
~~~~~\tilde{x}\rightarrow \hat{x}, ~\tilde{x}\in (\hat{x}+\tilde{\mathscr{S}})\cap \mathbb{R}^{n}_{>0}
\end{eqnarray}

In what follows, a special solution $G_e(\tilde {x})$ gratifying (\ref{D-PDE}) and (\ref{boundcon1}) is derived.

\begin{theorem}
For any CBP-CRNs ($\tilde{\mathcal{S}}, \tilde{\mathcal{C}}, \tilde{\mathcal{R}}, \tilde{\mathcal{K}}$) generated by a given complex balanced MAS ($\mathcal{S}, \mathcal{C}, \mathcal{R}, \mathcal{K}$), the following function
\begin{equation}\label{eq:ge}
G_e({\tilde{x}})=\sum^{n}_{j=1}d_j(\tilde{x}^*_j-\tilde{x}_j-\tilde{x}_j\ln {\frac{\tilde{x}^*_j}{\tilde{x}_j}})
\end{equation}
is a solution of the Lyapunov function PDEs (\ref{D-PDE}) and (\ref{boundcon1}).
\end{theorem}

\begin{pf}
$\forall ~\tilde{x} \in (\tilde{x}^*+ \tilde{\mathscr{S}})\cap{\mathbb{R}^{n}_{>0}}$, the gradient of the function (\ref{eq:ge}) is
\begin{eqnarray*}
\nabla G_e(\tilde{x})=D{\rm{Ln}}\frac{\tilde{x}}{\tilde{x}^*}
                     \xlongequal{\tilde{x}=D^{-1}x}D{\rm{Ln}}\frac{{x}}{x^*}
                     =D \nabla G(x).
\end{eqnarray*}
Substituting it into the left-hand side of (\ref{D-PDE}) and  (\ref{boundcon1}), and further utilizing conditions involved in the Definition \ref{df:CBP-CRNs}, we obtain
\begin{eqnarray}\label{CBP-CRNpf}
 &&\sum^{\tilde{r}}_{i=1}\tilde{k}_{i}\tilde{x}^{\tilde{v}_{\cdot
 i}}-
 \sum^{\tilde{r}}_{i=1}\tilde{k}_{i}\tilde{x}^{\tilde{v}_{\cdot
 i}}\exp\{(\tilde{v}'_{\cdot i}-\tilde{v}_{\cdot i})^\top \nabla G_e(\tilde{x})\} \notag \\
 &=&\sum^{r}_{i=1}k_i\prod^{n}_{j=1}d^{v_{ji}}_{j} (D^{-1}x)^{v_{\cdot i}}
 -\sum^{r}_{i=1}k_i\prod^{n}_{j=1}d^{v_{ji}}_{j} (D^{-1}x)^{v_{\cdot i}}\notag \\
&& \exp \{[ D^{-1}(v'_{\cdot i}-v_{\cdot i})]^\top D\nabla G\}\notag \\
 &=&\sum^{r}_{i=1}k_i x^{v_{\cdot i}}-\sum^{r}_{i=1}k_i x^{v_{\cdot i}}\exp\{(v'_{\cdot i}-v_{\cdot i})^\top\nabla G\} \notag \\
 &=&0.
\end{eqnarray}
As for the boundary condition (\ref{boundcon1}), it becomes

\begin{widetext}
\begin{eqnarray}\label{boundcon2}
&&\underset {\tilde{x}\rightarrow \hat{x}}{\lim}
\sum_{\{i|v_{\cdot i}\in \mathcal{C}_{\hat x}\}}\tilde{k}_i \tilde{x}^{v_{\cdot i}}
-\sum_{\{i|v'_{\cdot i}\in \mathcal{C}_{\hat x}\}}\tilde{k}_i \tilde{x}^{v_{\cdot i}}\exp \{(\tilde{v}'_{\cdot i}-\tilde{v}_{\cdot i})^\top\nabla {G}_e(\tilde{x})\}
\notag \\
&=&\underset {{x}\rightarrow \bar{x}'}{\lim}
\sum_{\{i|v_{\cdot i}\in \mathcal{C}_{\hat x}\}}
 k_i\prod^{n}_{j=1}d^{v_{ji}}_{j} (D^{-1}x)^{v_{\cdot i}}-\sum_{\{i|v'_{\cdot i}\in \mathcal{C}_{\hat x}\}}
 k_i\prod^{n}_{j=1}d^{v_{ji}}_{j}(D^{-1}x)^{v_{\cdot i}}
 \exp \{[D^{-1}(v'_{\cdot i}-v_{\cdot i)}]^\top D\nabla G\}\notag \notag \\
&=&\underset {{x}\rightarrow \bar{x}'}{\lim}
\sum_{\{i|v_{\cdot i}\in \mathcal{C}_{\hat x}\}}
k_ix^{v_{\cdot i}}
-\sum_{\{i|v'_{\cdot i}\in \mathcal{C}_{\hat x}\}}
k_ix^{v_{\cdot i}}\exp \{(v'_{\cdot i}-v_{\cdot i})^\top \nabla G\}
\notag \\
&=& 0,
~~~~\forall x\in (\bar{x}'+\mathscr{S})\cap \mathbb{R}^{n}_{>0}
\end{eqnarray}
\end{widetext}
Note that the last equality in (\ref{CBP-CRNpf}) and (\ref{boundcon2}) are true due to the fact that $G(x)$ is a solution of the Lyapunov function PDEs for complex balanced MASs shown in Theorem \ref{complexPDEth}.
$\Box$
\end{pf}

\begin{remark}
As stated in [\cite{Ke2019Complex}], a \emph{generalized pseudo-Helmholtz function} is defined as
\begin{equation*}
\tilde{G}(x)=\sum^{\tau_p}_{i=\tau_1} d_i(x^*_i-x_i-x_i\ln{\frac{x^*_i}{x^*_i}})
\end{equation*}
where $d_i\geq 0,~\tau_1\leq \cdots \leq \tau_p$ and $\{ \tau_1, \cdots, \tau_p \}\subseteq \mathcal{I}=\{1,\cdots, n\}$.

Observe that the expression of solution $G_e(\tilde{x})$ is a special case of $\tilde{G}(\tilde{x})$ while $i=1,\cdots, n$,  thus we call it the generalized pseudo-Helmholtz function in the remaining contents.
\end{remark}

Afterwards, we attend to show this solution $G_e(\tilde{x})$ can be taken as a Lyapunov function for the CBP-CRNs, and that CBP-CRNs inherit the properties on equilibria distribution from complex balanced MASs.
\begin{theorem}\label{CBP-CRNstability}
Given any CBP-CRN ($\tilde{\mathcal{S}}, \tilde{\mathcal{C}}, \tilde{\mathcal{R}}, \tilde{\mathcal{K}}$) with an equilibrium $\tilde{x}^* \in \mathbbold{R}^{n}_{>0}$ and generated from a complex balanced MAS ($\mathcal{S}, \mathcal{C}, \mathcal{R}, \mathcal{K}$), for any initial state $\tilde{x}_{0}\in\mathbb{R}_{>0}^{n}$, there exists a unique positive equilibrium $\tilde{x}^{+}$ of ($\tilde{\mathcal{S}}, \tilde{\mathcal{C}}, \tilde{\mathcal{R}}, \tilde{\mathcal{K}}$) such that $\tilde{x}^{+}\in(\tilde{x}_{0}+\tilde{\mathscr{S}})$, and the solution $G_e(\tilde{x})$ of (\ref{D-PDE}) and (\ref{boundcon1}) can act as a Lyapunov function to render $\tilde{x}^+$ to be locally asymptotically stable with respect to arbitrary initial condition in $(\tilde{x}_{0}+ \tilde{\mathscr{S}})\cap{\mathbb{R}^{n}_{>0}}$ near $\tilde{x}^+$.
\end{theorem}
\begin{pf}
 According to Proposition \ref{prop}, $x^{*}=D\tilde{x}^{*}$ is an equilibrium of ($\mathcal{S}, \mathcal{C}, \mathcal{R}, \mathcal{K}$). From Theorem \ref{thm:complex}, there exists a unique positive equilibrium $x^{+}\in({x}_{0}+\mathscr{S})$ of ($\mathcal{S}, \mathcal{C}, \mathcal{R}, \mathcal{K}$), which can lead to $\tilde{x}^{+}=D^{-1}x^{+}$ being a positive equilibrium in ($\tilde{ {x}}_{0}+\tilde{\mathscr{S}}$) of ($\tilde{\mathcal{S}}, \tilde{\mathcal{C}}, \tilde{\mathcal{R}}, \tilde{\mathcal{K}}$), and the uniqueness of $x^{+}$ guarantees that for $\tilde{x}^{+}$ based on Proposition \ref{prop}.\\
 Since the solution of Lyapunov PDEs is dissipative, we know $\dot {G_e}(\tilde{x})\leq 0$.
 Next we carry out the Hessian matrix of $G_e(\tilde{x})$, which
 is obtained as
\begin{equation*}
\begin{array}{lll}
\nabla^2 G_e(\tilde{x})&=&
\left(
\begin{array}{ccc}
d_1\tilde{x}^{-1}_1&& \\
&\ddots&\\
&&d_n\tilde{x}^{-1}_n\\
\end{array}
\right)>0,
\end{array}
\end{equation*}
Clearly, $\forall \tilde{x} \in (\tilde{x}_{0}+ \tilde{\mathscr{S}})\cap{\mathbb{R}^{n}_{>0}}$ near $\tilde{x}^+$, $G_e(\tilde{x})$ is twice differentiable and strictly convex, which means that $\dot {G_e}(\tilde{x})=0$ if and only if $\tilde{x}=\tilde{x}^+$.
Then the asymptotically stability follows immediately due to Theorem \ref{asymstability}.
$\Box$
\end{pf}

\begin{remark}
As stated in Remark \ref{rm:mark1}, there can be diverse networks induced by a single complex balanced system on account of various choices of producing matrices Ds, while for each of the former networks a Lyapunov candidate is naturally given. This means with above method, a lot more MASs other than complex balanced ones are equipped with a promising stability analysis approach. \\
At the same time, \cite{Fang2015Lyapunov}'s conjecture, that for any MAS with proper boundary complex set admitting a stable positive equilibrium, its corresponding Lyapunov function PDEs have a solution qualified as a Lyapunov function to suggest the system is locally asymptotically stable at the equilibrium, are confirmed for a class of CRNs with more general structures.
\end{remark}

\section{Illustration}
In this section an example is demonstrated to illustrate the above results.

Consider the following complex balanced network

~~~~~~~~~~~~~~~~
\xymatrix{
  2S_1 \ar[r]^{2}  &  3S_1 \ar[dl]^{1}    \\
  2S_1+S_2 \ar[u]^{2} &                    }
  \\
  \\
where
the species set =$\{S_1, ~S_2\}$, complex set=$\{v_{\cdot 1}=v'_{\cdot 3},~v'_{\cdot 1}=v_{\cdot 2}, ~v_{\cdot 3}=v'_{\cdot 2}\}$,
reaction set=$\{v_{\cdot 1}\rightarrow v'_{\cdot 1},~v_{\cdot 2}\rightarrow v'_{\cdot 2},~v_{\cdot 3}\rightarrow v'_{\cdot 3}\}$,

${v}_{\cdot 1}={v}'_{\cdot 3}=\left(
  \begin{array}{c}
    2 \\
    0 \\
  \end{array}
\right)$,
${v}'_{\cdot 1}={v}_{\cdot 2}=\left(
  \begin{array}{c}
    3 \\
    0 \\
  \end{array}
\right)$,
${v}'_{\cdot 2}={v}_{\cdot 3}=\left(
  \begin{array}{c}
    2 \\
    1 \\
  \end{array}
\right)$.
Then the dynamics are established as
\begin{eqnarray*}
~\left\{
\begin{array}{lll}
\dot x_1&&=2x^2_1-x^3_1 \\
\dot x_2&&=-x^3_1-2x^2_1x_2
\end{array}
\right.
\end{eqnarray*}
with an equilibrium $(x^*_1,x^*_2)=(2,1)$.

Setting $D={\rm diag} (\frac{1}{3},1)$, then according to the Definition \ref{df:CBP-CRNs} we have the following CBP-CRN
\begin{eqnarray}
2S_1 \stackrel{\frac{2}{9}}\longrightarrow 5S_1& \notag \\
3S_1 \stackrel{\frac{1}{27}}\longrightarrow S_2& \notag\\
2S_1+S_2 \stackrel{\frac{2}{9}}\longrightarrow 2S_1&
\end{eqnarray}
with the complexes being

$\tilde{v}_{\cdot 1 }=\tilde{v}'_{\cdot 3}=\left(
  \begin{array}{c}
    2 \\
    0 \\
  \end{array}
\right)$,
$\tilde{v}'_{\cdot 1}=\left(
  \begin{array}{c}
    5 \\
    0\\
  \end{array}
\right)$,
$\tilde{v}_{\cdot 2}=\left(
  \begin{array}{c}
    3 \\
    0 \\
  \end{array}
\right)$,
$\tilde{v}'_{\cdot 2}=\left(
  \begin{array}{c}
    0 \\
    1 \\
  \end{array}
\right)$,
$\tilde{v}_{\cdot 3}=\left(
  \begin{array}{c}
    2 \\
    1 \\
  \end{array}
\right)$.
Now the dynamical equations follow
\begin{eqnarray}\label{eg}
~\left\{
\begin{array}{lll}
\dot {\tilde{x}}_1&&=\frac{2}{3}\tilde{x}^2_1- \frac{1}{9}\tilde{x}^3_1 \\
\dot {\tilde{x}}_2&&=\frac{1}{27}\tilde{x}^3_1-\frac{2}{9}\tilde{x}^2_1\tilde{x}_2
\end{array}
\right.
\end{eqnarray}
and the equilibrium is computed as $(\tilde{x}^*_1, \tilde{x}^*_2)=(6,1)$, and its Lyapunov function PDEs are written as
\begin{eqnarray}\label{egpde}
&&\frac{2}{9}\tilde{x}^2_1\tilde{x}_2\left(1-\exp\{-\frac{\partial{f}}{\partial \tilde{x}_2}\}\right)
+\frac{2}{9}\tilde{x}^2_1 \left(1-\exp\{3\frac{\partial{f}}{\partial \tilde{x}_1}\}\right) \notag \\
&&+\frac{1}{27}\tilde{x}^3_1\left(1-\exp\{-3\frac{\partial{f}}{\partial \tilde{x}_1}+\frac{\partial{f}}{\partial \tilde{x}_2}\}\right)=0.
\end{eqnarray}
For simplicity, the boundary point is set as $(\hat{x}_1,\hat{x}_2)=(0,\hat{x}_2) $ with $\hat{x}_2>0$. Consequently, the boundary conditions are naturally satisfied.

Actually, it is not hard to clarify that
\begin{equation*}
f(\tilde{x}_1,\tilde{x}_2)=\frac{1}{3}(-\tilde{x}_1-\tilde{x}_1\ln {\frac{6}{\tilde{x}_1}})-\tilde{x}_2-\tilde{x}_2\ln {\frac{1}{\tilde{x}_2}}+3
\end{equation*}
is a solution of (\ref{egpde}).
Substituting
\begin{equation*}
\nabla f(\tilde{x}_1,\tilde{x}_2)=(\frac{1}{3}\ln\frac{\tilde{x}_1}{6}, \ln\tilde{x}_2)^\top
\end{equation*}
 into the left-hand side (L.H.S.) of (\ref{egpde}) yields
\begin{eqnarray*}
&&{\rm L.H.S. of}(\ref{egpde})\\
=&&\frac{2}{9}\tilde{x}^2_1 \tilde{x}_2-\frac{2}{9}\tilde{x}^2_1+\frac{2}{9}\tilde{x}^2_1-\frac{1}{27}\tilde{x}^3_1+\frac{1}{27}\tilde{x}^3_1-\frac{2}{9}\tilde{x}^2_1 \tilde{x}_2\\
=&&0.
\end{eqnarray*}
Additionally, Theorem \ref{CBP-CRNstability} tells us that $f(\tilde{x}_1,\tilde{x}_2)$ is valid for analyzing the asymptotic stability of $(6,1)$ (also see Fig.1).

\begin{figure}
\begin{center}
\includegraphics[width=8.4cm]{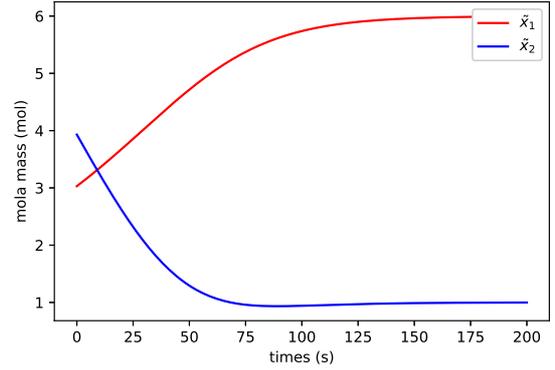}    
\caption{ state evolution with initial value (3,4)}
\end{center}
\end{figure}


\section{Conclusion}
This paper is devoted to seeking Lyapunov functions for the CBP-CRNs from the point view of establishing their Lyapunov function PDEs. The CBP-CRNs are a large class of networks generated from complex balanced systems and show some nice properties on equilibria distribution. 
Inspired by the work where the Lyapunov function PDEs for the complex balanced MASs are solved by the pseudo-Helmholtz  function $G(x)$,  we construct a special solution  $G_e(\tilde{x})$ of the developed Lyapunov function PDEs corresponding to the CBP-CRNs. Besides, it demonstrates that $G_e(\tilde{x})$ plays the role as a Lyapunov function for the CBP-CRNs to be asymptotically stable at the equilibrium. \\
Apparently, these results have extended \cite{Fang2015Lyapunov}'s work of Lyapunov function PDEs to reaction networks with more general structures, and can be taken as an argument to their conjecture that Lyapunov function PDEs can serve for any CRNs. Meanwhile, the performance of Lyapunov functions on CRNs with other structures still needs to be explored. Our future work will concentrate on such systems as the linear reaction networks and autocatalytic reaction networks.
\\~


\normalem








\end{document}